# *H*∞ Consensus of nonlinear complex multi-agent systems using dynamic output feedback controller: An LMI approach

Amir Amini, Ali Azarbahram, Mahdi Sojoodi[*]

*Advanced Control Systems Laboratory, School of Electrical and Computer Engineering, Tarbiat Modares University, Tehran, Iran. (\*Corresponding Author)*

**Abstract:**
This paper investigates a new method for consensus in a group of nonlinear complex multi-agent systems using fixed-order non-fragile dynamic output feedback controller, via an LMI approach. The proposed scheme is decentralized in the sense that each agent relies on the relative output information among the adjacent agents. The consensus based controllers are designed to minimize the effects of nonlinear terms of the agents as well as external disturbances. Converting consensus problem to stabilization of an equivalent augmented system using proper transformations, Lyapunov stability theorem is applied to obtain unknown controller parameters in order to guarantee consensus and simultaneously acquire considered control objectives. Finally, to demonstrate the effectiveness of the proposed algorithm and compare with similar earlier researches, a numerical example on a multi-agent system consisting of single link flexible manipulators is carried out.

## 1. Introduction

Coordination and cooperation of complex multi-agent systems (MASs) which is composed of multiple autonomous agents with mutual interactions has been an interesting research topic in the last decade. Local communication between individual agents in the network may result in certain desirable global behaviors such as formation control, flocking, leader following and consensus. To reach an agreement regarding the states of all agents in a networked multi-agent system is defined as consensus. Recently, consensus has been widely studied because of its noteworthy applications such as cooperative control of unmanned airplanes and underwater vehicles, flocking of mobile robots, sensor networks and rendezvous [1]. The consensus problem has been investigated from different viewpoints including various agents' dynamics, constraints on network topologies, control strategies and design methods. Many aspects of consensus algorithms is analyzed for multi-agent networked systems in [2], including convergence methods, robustness to topology changes and communication delays. In [3] an adaptive distributed consensus problem of continuous time multi-agent systems with general linear dynamics for both cases of leaderless and leader following consensus with time varying communication topology is designed. Consensus in a network of agents with linear or linearized general dynamics on a consensus region using static state feedback is introduced in [4]. In [5] a consensus protocol in a multi-agent system composed of identical agents with general linear dynamics of any order is presented considering static and switching network topology. Sufficient conditions for consensus in linear second order multi-agent systems with directed network topologies are derived in [6]. Conditions on consensus algorithms in multi-agent systems with second order dynamics containing nonlinear terms under directed topologies have been studied in [7]. In [8] leader following consensus problem of second-order nonlinear multi-agent systems is investigated via pinning control algorithm. Using an adaptive distributed protocol considering neighboring states, consensus problem for general linear multi-agent systems and nonlinear agents satisfying Lipschitz condition is addressed in [9]. In [10] consensus problem for low order multi-agent system with inherent Lipschitz nonlinear dynamics under directed communication topology is considered. A delay-dependent distributed method is studied in [11] to insure consensus in nonlinear second order multi-agent systems under switching topology. Nevertheless, a systematic procedure which ensures consensus for multi-agent systems with general linear dynamics and a nonlinear term which satisfies sufficient condition is essential. A static feedback consensus protocol for identical networked systems with linear high-order dynamics is introduce in [12] based on only output measurements of the neighboring agents. State consensus problem of linear multi-agent systems by simultaneous stabilization techniques is proposed in [13], using neighbor-based static and dynamic output-feedback control algorithm. Another static output feedback controller for $H_\infty$ consensus in a class of linear time-varying multi-agent system with both missing information and model uncertainties is designed in [14].

However, in most practical cases the linearity of multi-agent system dynamics is contradicted, hence methods which guarantee consensus in multi-agent systems considering nonlinear dynamics are required. Consensus tracking in second-order multi-agent systems with bounded disturbances using a nonlinear distributed scheme is studied in [15]. The robust consensus tracking problem in [16] concerns multi-agent systems with uncertain Lur'e-type nonlinear dynamics with a fixed network topology utilizing a protocol based on relative state information among the neighboring sets. However, for seeking consensus in multi-agent systems, a general linear dynamic controller which uses only measurable system states is preferable.

Utilizing matrix inequality techniques, in [17] a consensus algorithm for general linear multi-agent systems using local state information is proposed. The approach is based on graph Laplacian of the network and reduces the consensus problem to stabilization of the augmented system. This algorithm is expanded to design dynamic output feedback controller in [18], however, the final matrix inequalities are bilinear, therefore, iterative and inaccurate methods of solving Bilinear Matrix Inequalities (BMIs) are needed. Moreover, [19] investigates the output consensus problem of high-order multi-agent systems with external disturbances and proposes a distributed protocol using the neighbors' output information. Nonetheless, some predefined conditions on controller design procedure are required to convert the initially obtained BMI to an LMI form.

To the best of our knowledge, no result is available on the LMI-based robust fixed-order decentralized dynamic output feedback consensus algorithm design for multi-agent systems with nonlinear dynamics. Furthermore, the existence of uncertainties in controller implementation has also been taken into account to cover controllers' fragility. The method, offers the maximum possible bounds on the agents' nonlinearities and simultaneously attenuate disturbance with a prescribed attenuation level, via LMI framework.

The rest of this paper is organized as follows: Section 2 provides statement and some preliminary backgrounds. In Section 3, the problem is formulated, the proposed controller design method is presented and the theoretical results are provided in the form of a theorem and a corollary. The case studies are discussed in Section 4. Finally, some conclusion remarks are given in Section 5.

## 2. Statement and preliminaries

Some mathematical notations are defined here to be used throughout this paper. We use $\|.\|$ to define the Euclidean norm. Let $\mathbf{1}_n$ and $\mathbf{0}_n$ denote $n \times 1$ column vectors with all elements to be ones and zeros, respectively. Moreover, $I_n$ indicates the $n \times n$ identity matrix while $\mathbf{J}_n$ is all-ones matrix of $n \times n$ dimension. For the kronecker product the symbol $\otimes$ is used and for the two matrices of $A$ and $B$ with the same dimensions $m \times n$, the Hadamard product, $A \circ B$, is a matrix, of the same dimension as the operands, with elements given by:

$(A \circ B)_{ij} = (A)_{ij}.(B)_{ij}$

A symmetric positive definite matrix $M$ is defined by $M > 0$ on the other hand, pseudo inverse of a given non-square matrix $A_{n \times m}$ is shown by $A^\uparrow$.

Graph theory is widely used to model the communication topology of the network of multi-agent systems. Therefore, some necessary and required concepts regarding graph theory are introduced. The communication network of a multi-agent system composed of $N$ agents is usually modeled using a directed graph $G = \{V, E, A\}$, where $V = \{1, 2, ..., N\}$ is a set of $N$ integers in which the $i-th$ vertex indicates the $i-th$ agent while the edge set is $E \subseteq V \times V$. If the information of the $j-th$ agent is available for the $i-th$ one, then the pair $(j, i)$ is an element of $E$ and is depicted by $j \to i$ in graph representation. accordingly, $(j, i) \in E$ is equivalent to $(i, j) \in E$ .in an undirected graph.

The set of neighbors of agent $i$ is denoted by $N_i = \{j \in V : (i, j) \in E\}$ and $A = [a_{ij}] \in \Re^{n \times n}$ is a weighted adjacency matrix of $G$, where $a_{ii} = 0$ and $a_{ij} = 1$ if $(i, j) \in V$ or $a_{ij} = 0$ otherwise. The Laplacian of $G$ is defined as $L = D - A$, where:

$$D = diag(\deg_1, ..., \deg_N), \quad \deg_i = \sum_{j=1}^{N} a_{ij}$$

In other words the Laplacian graph is:

$$l_{ii} = \sum_{j \in N_i} a_{ij}, \quad l_{ij} = -a_{ij} \quad \forall i \neq j, \quad i, j = 1, 2, ..., N \tag{1}$$

A sequence of edges in the form of $(i_1, i_2), (i_2, i_3), ..., (i_{N-1}, i_N)$ with $i_K \in V$, creates a path in an undirected graph. If there exists a path between each two vertices, the undirected graph is connected. It can be proved that zero is always a simple eigenvalue of the Laplacian matrix of a connected graph and $\mathbf{1}_n$ is its corresponding eigenvector while, all nonzero eigenvalues are positive.

## 3. Problem formulations

Consider a multi-agent system comprises $N$ agents. The nonlinear dynamics of the $i-th$ agent can be represented by:

$$\begin{aligned}\dot{x}_i &= \bar{A} x_i + B_i u_i + h_i(t, x_i) + \xi_i(t) \\ y_i &= \bar{C} x_i \quad\quad\quad i = 1, ..., N\end{aligned} \tag{2}$$

where $x_i \in \Re^n$, $u_i \in \Re^m$, $y_i \in \Re^q$ and $\xi_i(t)$ are agent states, control input, output and the external disturbance of the $i-th$ agent, while $\bar{A}, B_i$ and $\bar{C}$ are known constant matrices of appropriate dimensions. The pairs of $(\bar{A}, B_i)$ and $(\bar{A}, \bar{C})$ are assumed to be controllable and observable, respectively. Then, $h_i(t, x_i)$ contains all possible nonlinearities and uncertainties of the $i-th$ agent which is a piecewise continuous vector function satisfies the following quadratic inequality in its domain of continuity [20]:

$$h_i(t, x_i)^T h_i(t, x_i) \leq \bar{\alpha}_i^2 x_i^T \bar{H}_i^T \bar{H}_i x_i \quad i = 1, ..., N \tag{3}$$

where $\bar{\alpha}_i > 0$ is a bounding parameter and supposed to be maximized while $\bar{H}_i$ is a constant $v_i \times n$ matrix.

Thus the entire multi-agent system can be represented as:

$$\begin{aligned}\dot{x} &= Ax + Bu + h(t, x) + \xi(t) \\ y &= Cx\end{aligned} \tag{4}$$

where $x = [x_1^T \ x_2^T \ ... \ x_N^T]^T$, $u = [u_1^T \ u_2^T \ ... \ u_N^T]^T$ and $y = [y_1^T \ y_2^T \ ... \ y_N^T]^T$ are the global states, control inputs and outputs, respectively. $A = I_N \otimes \bar{A}$ is the global system matrix, $B = \text{diag}(B_1,...,B_N)$ is the global input matrix and $C = I_N \otimes \bar{C}$ is the global output matrix and $h(t,x) = [h_1(t,x_1)^T \ \cdots \ h_N(t,x_N)^T]^T$ is the global nonlinear vector function. Let $\bar{H} = [\bar{H}_1^T \ \cdots \ \bar{H}_N^T]^T$ where $\bar{H}_i$, $i = 1,...,N$ are defined in (3), and $\bar{\Gamma}_1 = \text{diag}\{\bar{\gamma}_{11}I_{v_1},...,\bar{\gamma}_{1N}I_{v_N}\}$ with $\bar{\gamma}_{1i} = \bar{\alpha}_i^{-2}$. Then, it is always possible to find matrices $H$ and $\Gamma_1$ such that [20]:

$$h(t,x)^T h(t,x) \leq x^T \bar{H}^T \bar{\Gamma}_N^{-1} \bar{H} x \leq x^T H^T \Gamma_N^{-1} H x \tag{5}$$

where $H = \text{diag}(H_1, H_2, \cdots, H_N)$, with $v_i \times n_i$ matrices of $H_i$, and:

$$\Gamma_N = \text{diag}\{\gamma_{11}I_{v_1},...,\gamma_{1N}I_{v_N}\}, \ \gamma_{1i} > 0, \ i = 1,...,N. \tag{6}$$

The sufficient condition for (5) gives the following condition for matrices $H$ and $\Gamma_N$ [6]:

$$\lambda_{\max}(\bar{H}^T \bar{H}) \min_i \bar{\gamma}_{1i} \leq \max_i \gamma_{1i} \min_i \lambda_{\min}(H_i^T H_i) \tag{7}$$

The objective is to design a protocol such that the multi-agent system (2) asymptotically reaches global consensus on the states, which is expressed by:

$$\lim_{t \to \infty} |x_i(t) - x_j(t)| = 0, \quad \forall i \neq j \tag{8}$$

For this purpose, a decentralized dynamic output feedback controller is proposed based on the adjacent output information between the agents. Hence, the designed controllers may be sensitive to errors in the control parameters, the applied controllers usually encounter difficulties regarding inaccuracies. This issue compels the designers to adopt a procedure in order to avoid possible undesirable or unstable closed loop performance at implement stage. Therefore, the following decentralized non-fragile dynamic output feedback controller is presented to guarantee consensus among agents in multi-agent system.

$$\begin{aligned} \dot{x}_C &= (A_C + \Delta A_C)x_C + (B_C + \Delta B_C)L_q y \\ u &= (C_C + \Delta C_C)x_C + (D_C + \Delta D_C)L_q y \end{aligned} \tag{9}$$

where $x_C \in \Re^{Nn_C} = [x_{C_1}^T \ ... \ x_{C_N}^T]^T$, $L_q = L \otimes I_q$. Additionally $A_C$, $B_C$, $C_C$ and $D_C$ are the global controller unknown parameters of the form:

$$\begin{aligned} A_C &= \text{diag}(A_{C_1},...,A_{C_N}) \\ B_C &= \text{diag}(B_{C_1},...,B_{C_N}) \\ C_C &= \text{diag}(C_{C_1},...,C_{C_N}) \\ D_C &= \text{diag}(D_{C_1},...,D_{C_N}) \end{aligned} \tag{10}$$

The perturbations uncertainties in controller parameters satisfy the norm bounded conditions of the form:

$$\|\Delta A_C\| \leq \delta_{A_C}, \ \|\Delta B_C\| \leq \delta_{B_C}, \ \|\Delta C_C\| \leq \delta_{C_C}, \ \|\Delta D_C\| \leq \delta_{D_C} \tag{11}$$

The closed loop augmented multi-agent system composed of equations (4) and (9) can be, therefore, represented by:

$$\begin{bmatrix} \dot{x} \\ \dot{x}_C \end{bmatrix} = \begin{bmatrix} A + B(D_C + \Delta D_C)L_q C & B(C_C + \Delta C_C) \\ (B_C + \Delta B_C)L_q C & A_C + \Delta A_C \end{bmatrix} \begin{bmatrix} x \\ x_C \end{bmatrix} + \begin{bmatrix} h(t,x) \\ 0 \end{bmatrix} + \begin{bmatrix} \xi(t) \\ 0 \end{bmatrix} \quad (12)$$

Since our idea is to reduce consensus problem in (12) to a stabilization problem, it is necessary to state two essential obstacles toward the controller parameters attainment and consensus assurance. First, if a controller is designed for (12) such that the overall closed loop system is stable then the initial states only converge to zero, however compelling corresponding states to the same value as intended in consensus problem, it is necessary to define a transformation as:

$$\tilde{x} = L_n x, \quad L_n = L \otimes I_n \quad (13)$$

The interference in (13) is based on the fact that $\tilde{x}$ provides relative state difference of the multi-agent system according to graph Laplacian matrix, hence if $\tilde{x}$ converges to zero, then the consensus definition in (8) is achieved among the agents. Nonetheless, as the second problem, Lyapunov theorem, which is the major criterion for stability analysis, requires a full rank matrix system in state-space representation that is not satisfied in (12) after applying the transformation (13), due to rank deficiency of $L_n$. Therefore, removing one of the rows of $L$ will eliminate system redundancy accordingly a reduced full rank system matrix is achieved by modifying (13):

$$x_r = \hat{L}_n x, \quad \hat{L}_n = \hat{L} \otimes I_n \quad (14)$$

where $\hat{L} \in \Re^{(N-1) \times N}$ is the independent rows of $L$. In order to obtain state space representation with respect to $x_r$, the following two lemmas are introduced based on graph Laplacian properties.

Lemma 1 [15]: $L_q C = C L_n$

*Proof*: $L_q C = (L \otimes I_q)(I_N \otimes \overline{C}) = (LI_N) \otimes (I_q \overline{C}) = (I_N L) \otimes (\overline{C} I_n) = (I_N \otimes \overline{C})(L \otimes I_n) = C L_n$

Lemma 2 [15]: $L_n A = A L_n$

*Proof*: $L_n A = (L \otimes I_n)(I_N \otimes \overline{A}) = (LI_N) \otimes (I_n \overline{A}) = (I_N L) \otimes (\overline{A} I_n) = (I_N \otimes \overline{A})(L \otimes I_n) = A L_n$

Using $\hat{L}_n$ and subsequently $x_r \in \Re^{(N-1) \times n}$, the only remaining challenge is to match matrices dimensions in the resulting system matrix. Therefore, $A_r = I_{N-1} \otimes \overline{A}$ and $C_r$ will replace $A$ and $C$, respectively, where $C_r = CL_n \hat{L}_n^{\uparrow}$ hence the alteration of $C$ to $C_r$ is necessary to satisfy $C_r \hat{L}_n x = CL_n x$. Noticing all the above mentioned considerations, the closed loop system of (12) transforms to (15) after replacing $x_r$, $A_r$, $C_r$ and applying Lemma 1 and Lemma 2:

$$\dot{x}_{Cl} = A_{Cl} x_{Cl} + h_r + \xi_r \quad (15)$$

where

$$x_{Cl} = \begin{bmatrix} x_r \\ x_C \end{bmatrix}, \quad A_{Cl} = \begin{bmatrix} A_r + \hat{L}_n B(D_C + \Delta D_C)C_r & \hat{L}_n B(C_C + \Delta C_C) \\ (B_C + \Delta B_C)C_r & A_C + \Delta A_C \end{bmatrix}, \quad h_r = \begin{bmatrix} \hat{L}_n h(t,x) \\ 0 \end{bmatrix}, \quad \xi_r = \begin{bmatrix} \hat{L}_n \xi(t) \\ 0 \end{bmatrix}.$$

Separating $A_{Cl}$ into two portions of certain and uncertain parameters as:

$$A_{Cl} = A_\Delta + A_\Phi \quad (16)$$

where

$$A_\Phi = \begin{bmatrix} A_r + \hat{L}_n BD_C C_r & \hat{L}_n BC_C \\ B_C C_r & A_C \end{bmatrix}, \quad A_\Delta = \begin{bmatrix} \hat{L}_n B\Delta D_C C_r & \hat{L}_n B\Delta C_C \\ \Delta B_C C_r & \Delta A_C \end{bmatrix},$$

simplifies our procedure towards stating a theorem which gives the sufficient condition on designing a fixed-order dynamic output feedback controller for consensus of a nonlinear multi-agent system ensuring $H_\infty$ performance of the closed-loop system (15).

**Definition 1.** $H_\alpha$ **class**: For any given $\alpha = [\alpha_1^T \cdots \alpha_N^T]^T$ and matrix $H$ with appropriate dimension, $H_\alpha$ is a class of piecewise-continues functions defined by:

$$H_\alpha = \{h(t,x) | h \in \Re^n, h^T h \leq x^T H^T \alpha^T \alpha H x \quad \text{in the domain of continuity}\} \tag{17}$$

Therefore, if $h \in H_\alpha$, then $h(t,0) = 0$, which implies that $x = 0$ is the equilibrium of system (15).

**Definition 2.** Robust stability with vector degree $\alpha = [\alpha_1^T \cdots \alpha_N^T]^T$: System (15) with bounded nonlinear uncertainty given in (5) is called robustly stable with vector degree $\alpha$ if the equilibrium $x = 0$ is globally asymptotically stable for all $h(t,x) \in H_\alpha$.

Definitions 1 and 2 are the vector expansion of those in [20] and are supposed to be utilized here.

**Theorem 1**: Considering augmented closed-loop system in (15) and given matrices $R$ and $Q$ described by $R = diag\,(r_1,...,r_{N-1})$ and $Q = diag\,(q_1,...,q_N)$ with $q_i$ and $r_i$ as prescribed symmetric positive definite weighting matrices, if positive definite matrix $P = P^T$ in the form of

$$P = diag\,(P_S, P_C), \quad P_S = I_{N-1} \otimes \bar{p}_s, \quad P_C = diag\,(P_{C_1},...,P_{C_N}), \tag{18}$$

and matrices $Z = diag\,(z_1,...z_N)$, $W = diag\,(w_1,...w_N)$, $f_i$ and $k_i$ for $i = 1,..,N$ and positive scalars $\tau_j$ for $j = 1,..,5$ exist such that the following minimization problem becomes feasible:

$$\text{minmize} \quad \rho^2 + \sum_{i=1}^{N-1} \hat{\gamma}_i$$

$$\text{subject to} \quad \begin{bmatrix} \Phi + \Delta + \tilde{Q} & P & P\begin{bmatrix}\hat{L}_n B\\0\end{bmatrix} & P\begin{bmatrix}\hat{L}_n B\\0\end{bmatrix} & P\begin{bmatrix}0\\I\end{bmatrix} & P\begin{bmatrix}0\\I\end{bmatrix} & P\begin{bmatrix}I\\0\end{bmatrix} & \tau_5\begin{bmatrix}\hat{H}\\0\end{bmatrix} \\ * & -\rho^2 I & 0 & 0 & 0 & 0 & 0 & 0 \\ * & * & -\tau_1 I & 0 & 0 & 0 & 0 & 0 \\ * & * & * & -\tau_2 I & 0 & 0 & 0 & 0 \\ * & * & * & * & -\tau_3 I & 0 & 0 & 0 \\ * & * & * & * & * & -\tau_4 I & 0 & 0 \\ * & * & * & * & * & * & -\tau_5 I & 0 \\ * & * & * & * & * & * & * & -\hat{\Gamma} \end{bmatrix} < 0 \tag{19}$$

where

$$\Phi = \begin{bmatrix} A_r^T P_S + P_S A_r + \Pi_D C_r + C_r^T \Pi_D^T & C_r^T Z^T + \Pi_C \\ * & W + W^T \end{bmatrix},$$

$$\Delta = \begin{bmatrix} \|C_r\|(\tau_2 \delta_{D_C}^2 + \tau_4 \delta_{B_C}^2) I_{n(N-1)} & 0 \\ * & (\tau_1 \delta_{C_C}^2 + \tau_3 \delta_{A_C}^2) I_{n_C N} \end{bmatrix},$$

$$\tilde{Q} = diag(R_{n(N-1)}, Q_{n_C N}), \quad \hat{H} = diag(H_1, ..., H_{N-1}), \quad \hat{\Gamma} = diag(\hat{\gamma}_1 I_{v_1}, ..., \hat{\gamma}_{N-1} I_{v_{N-1}}),$$

(20)

with

$$\Pi_D = (\hat{L} \otimes \mathbf{1}_n) \circ (\mathbf{1}_{N-1} \otimes [k_1 \ ... \ k_N]),$$
$$\Pi_C = (\hat{L} \otimes \mathbf{J}_{n \times nc}) \circ (\mathbf{1}_{N-1} \otimes [f_1 \ ... \ f_N]),$$

(21)

then, the non-fragile dynamic output feedback controller parameters of:

$$A_C = P_C^{-1} W,$$
$$B_C = P_C^{-1} Z,$$
$$C_{C_1} = B_1^{\uparrow} \bar{p}_s^{-1} f_1, \ C_{C_2} = B_2^{\uparrow} \bar{p}_s^{-1} f_2, \ ..., \ C_{C_N} = B_N^{\uparrow} \bar{p}_s^{-1} f_N,$$
$$D_{C_1} = B_1^{\uparrow} \bar{p}_s^{-1} k_1, \ D_{C_2} = B_2^{\uparrow} \bar{p}_s^{-1} k_2, \ ..., \ D_{C_N} = B_N^{\uparrow} \bar{p}_s^{-1} k_N,$$

(22)

make the overall closed-loop system (15) stable with the robustness degree vector $[1/\sqrt{\gamma_1} \ \cdots \ 1/\sqrt{\gamma_{N-1}}]^T$ and simultaneously guarantee $H_\infty$ performance.

***Proof***: The idea adopted to poof Theorem 1 concerns the robust stabilization analysis of the system (15) within an optimization framework. Considering $H_\infty$ disturbance attenuation performance index with respect to $\xi_r$, the subsequent inequality must be satisfied [21]:

$$\dot{V}(x_{Cl}) + x_{Cl}^T \tilde{Q} x_{Cl} - \rho^2 \xi_r^T \xi_r \leq 0 \tag{23}$$

where $V(x_{Cl})$ is a positive quadratic Lyapunov candidate function and $\tilde{Q}$, as defined in (20), weights both the states of multi-agent system and the controllers using predefined $R$ and $Q$ matrices. The Lyapunov quadratic function candidate has the form of:

$$V(x_{Cl}) = x_{Cl}^T P x_{Cl}, \quad P > 0, \tag{24}$$

Then according to (15), $\dot{V}(x_{Cl})$ can be expanded as:

$$\dot{V}(x_{Cl}) = \dot{x}_{Cl}^T P x_{Cl} + x_{Cl}^T P \dot{x}_{Cl}$$
$$= (A_{Cl} x_{Cl} + h_r + \xi_r)^T P x_{Cl} + x_{Cl}^T P (A_{Cl} x_{Cl} + h_r + \xi_r). \tag{25}$$

With decomposition of $A_{Cl}$ introduced in (16), the equation (25) can be rewritten as:

$$\dot{V}(x_{Cl}) = \left( A_\Phi x_{Cl} + \begin{bmatrix} \hat{L}_n B(z_1 + z_2) + z_5 \\ z_3 + z_4 \end{bmatrix} \right)^T P x_{Cl} + x_{Cl}^T P \left( A_\Phi x_{Cl} + \begin{bmatrix} \hat{L}_n B(z_1 + z_2) + z_5 \\ z_3 + z_4 \end{bmatrix} \right)$$
$$+ \xi_r^T P x_{Cl} + x_{Cl}^T P \xi_r \tag{26}$$

where $z_1 = \Delta C_C x_C$, $z_2 = \Delta D_C C_r x_r$, $z_3 = \Delta A_C x_C$, $z_4 = \Delta B_C C_r x_r$ and $z_5 = \hat{L}_n h(t,x)$. Introducing $\Psi = [x_{Cl}^T \ \xi_r^T \ z_1^T \ z_2^T \ z_3^T \ z_4^T \ z_5^T]^T$, the equation in (26) is rearranged:

$$\dot{V}(x_{Cl}) = \Psi^T \begin{bmatrix} A_\Phi^T P + PA_\Phi & P & P\begin{bmatrix}\hat{L}_n B\\0\end{bmatrix} & P\begin{bmatrix}\hat{L}_n B\\0\end{bmatrix} & P\begin{bmatrix}0\\I\end{bmatrix} & P\begin{bmatrix}0\\I\end{bmatrix} & P\begin{bmatrix}I\\0\end{bmatrix} \\ * & 0 & 0 & 0 & 0 & 0 & 0 \\ * & * & 0 & 0 & 0 & 0 & 0 \\ * & * & * & 0 & 0 & 0 & 0 \\ * & * & * & * & 0 & 0 & 0 \\ * & * & * & * & * & 0 & 0 \\ * & * & * & * & * & * & 0 \end{bmatrix} \Psi \qquad (27)$$

Substituting (27) into (23), the subsequent matrix inequality is obtained:

$$\Psi^T \begin{bmatrix} A_\Phi^T P + PA_\Phi + \tilde{Q} & P & P\begin{bmatrix}\hat{L}_n B\\0\end{bmatrix} & P\begin{bmatrix}\hat{L}_n B\\0\end{bmatrix} & P\begin{bmatrix}0\\I\end{bmatrix} & P\begin{bmatrix}0\\I\end{bmatrix} & P\begin{bmatrix}I\\0\end{bmatrix} \\ * & -\rho^2 I & 0 & 0 & 0 & 0 & 0 \\ * & * & 0 & 0 & 0 & 0 & 0 \\ * & * & * & 0 & 0 & 0 & 0 \\ * & * & * & * & 0 & 0 & 0 \\ * & * & * & * & * & 0 & 0 \\ * & * & * & * & * & * & 0 \end{bmatrix} \Psi < 0 \qquad (28)$$

According to (11), the following norm bounds inequalities with respect to $z_1$, $z_2$, $z_3$ and $z_4$ are supposed to be satisfied:

$$z_1^T z_1 \leq \delta_{C_C}^2 x_C^T x_C, \qquad z_2^T z_2 \leq \|C_r\|\delta_{D_C}^2 x_r^T x_r,$$
$$z_3^T z_3 \leq \delta_{A_C}^2 x_C^T x_C, \qquad z_4^T z_4 \leq \|C_r\|\delta_{B_C}^2 x_r^T x_r. \qquad (29)$$

Furthermore, according to (5) one can conclude:

$$z_5^T z_5 \leq x_r^T \hat{H}^T \Gamma_{N-1}^{-1} \hat{H} x_r. \qquad (30)$$

Then, the combination of equations (29) and (30) with respect to $\Psi$ is:

$$\Psi^T \begin{bmatrix} \begin{bmatrix} \|C_r\|(\delta_{D_C}^2 + \delta_{B_C}^2)I & 0 \\ * & (\delta_{C_C}^2 + \delta_{A_C}^2)I \end{bmatrix} + \begin{bmatrix} \hat{H}^T \Gamma_{N-1}^{-1} \hat{H} & 0 \\ 0 & 0 \end{bmatrix} & 0 & 0 & 0 & 0 & 0 & 0 \\ * & 0 & 0 & 0 & 0 & 0 & 0 \\ * & * & -I & 0 & 0 & 0 & 0 \\ * & * & * & -I & 0 & 0 & 0 \\ * & * & * & * & -I & 0 & 0 \\ * & * & * & * & * & -I & 0 \\ * & * & * & * & * & * & -I \end{bmatrix} \Psi \geq 0 \qquad (31)$$

Repeatedly using *S-procedure* ([22]) for (28) and (31) result in:

$$\Psi^T \Omega \Psi < 0 \qquad (32)$$

where:

$$\Omega = \begin{bmatrix} A_\Phi^T P + P A_\Phi + \tilde{Q} + \Delta + \tau_5 \begin{bmatrix} \hat{H}^T \Gamma_{N-1}^{-1} \hat{H} & 0 \\ 0 & 0 \end{bmatrix} & P & P\begin{bmatrix} \hat{L}_n B \\ 0 \end{bmatrix} & P\begin{bmatrix} \hat{L}_n B \\ 0 \end{bmatrix} & P\begin{bmatrix} 0 \\ I \end{bmatrix} & P\begin{bmatrix} 0 \\ I \end{bmatrix} & P\begin{bmatrix} I \\ 0 \end{bmatrix} \\ * & -\rho^2 I & 0 & 0 & 0 & 0 & 0 \\ * & * & -\tau_1 I & 0 & 0 & 0 & 0 \\ * & * & * & -\tau_2 I & 0 & 0 & 0 \\ * & * & * & * & -\tau_3 I & 0 & 0 \\ * & * & * & * & * & -\tau_4 I & 0 \\ * & * & * & * & * & * & -\tau_5 I \end{bmatrix} \quad (33)$$

Thus, if $\Omega < 0$, the condition in (32) is satisfied. However, the matrix equation (33) is not an LMI due to several multiplication of variables. To deal with this problem, we first expand $A_\Phi^T P + P A_\Phi$ according to $P$ defined in (18).

$$A_\Phi^T P + P A_\Phi = \begin{bmatrix} A_r^T P_s + P_s A_r + C_r^T D_C^T B^T \hat{L}_n P_s + P_s \hat{L}_n B D_C C_r & P_s \hat{L}_n B C_C + C_r^T B_C^T P_C \\ * & P_C A_C + A_C^T P_C \end{bmatrix}$$

where

$$P_s \hat{L}_n B D_C = (\hat{L} \otimes \mathbf{1}_n) \circ (\mathbf{1}_{N-1} \otimes [\bar{p}_s B_1 D_{C1}, ..., \bar{p}_s B_N D_{CN}]),$$
$$P_s \hat{L}_n B C_C = (\hat{L} \otimes \mathbf{J}_{n \times nc}) \circ (\mathbf{1}_{N-1} \otimes [\bar{p}_s B_1 C_{C1}, ..., \bar{p}_s B_N C_{CN}]).$$

Now by changing variables as:

$W = P_C A_C$, $Z = P_C B_C$, $k_i = \bar{p}_s B_i D_{Ci}$ and $f_i = \bar{p}_s B_i C_{Ci}$

$\Pi_D$ and $\Pi_C$ in (21), and subsequently $\Phi$ in (20), are obtained respectively. Finally a schur complement is applied and introducing new variable $\hat{\Gamma} = \tau_5 \Gamma_{N-1}$ completes the proof.□

**Corollary 1**: Consider open-loop multi-agent system (4) and the controller (9) without uncertainties, i.e., $\Delta A_c = 0$, $\Delta B_c = 0$, $\Delta C_c = 0$ and $\Delta D_c = 0$. Given matrices $R$ and $Q$ described in Theorem 1, if positive definite matrix $P = P^T$ of the form (18) and matrices $Z = diag(z_1, ... z_N)$, $W = diag(w_1, ... w_N)$, $f_i$ and $k_i$ for $i = 1, .., N$ and positive scalar $\tau_1$, exist such that the following minimization problem is feasible:

$$\text{minmize} \quad \rho^2 + \sum_{i=1}^{N-1} \hat{\gamma}_i$$

$$\text{subject to} \begin{bmatrix} \Phi + \tilde{Q} & P & P\begin{bmatrix} I \\ 0 \end{bmatrix} & \tau_1 \begin{bmatrix} \hat{H} \\ 0 \end{bmatrix} \\ * & -\rho^2 I & 0 & 0 \\ * & * & -\tau_1 I & 0 \\ * & * & * & -\hat{\Gamma} \end{bmatrix} < 0 \quad (34)$$

where $\Phi$, $\tilde{Q}$, $\hat{H}$ and $\hat{\Gamma}$ are all defined in (20), then, the output feedback controller parameters given in (22), make the overall closed-loop system stable with the robustness degree vector $[1/\sqrt{\gamma_1} \quad \cdots \quad 1/\sqrt{\gamma_{N-1}}]^T$ and simultaneously guarantee $H_\infty$ performance.

***Proof:*** The proof procedure is similar to that of Theorem 1, Assuming $\Psi = [x_{Cl}^T \ \xi_r^T \ z_5^T]^T$.

**Remark 1:** In the case that the order of dynamic output feedback controller in (9) is set to be zero, then Corollary 1 and Theorem 1 are reduced to static output feedback controller design. In addition, if all the states are available, assuming $\bar{C} = I_n$, a dynamic static feedback controller can be designed by the procedure in Corollary 1 and Theorem 1.

## 4. Simulation

In this section, to validate the effectiveness of the proposed scheme, the designed dynamic output feedback controller in Theorem 1 and Corollary 1 are applied for a multi-agent system consisting single arm link manipulators. Using Matlab/Yalmip package for solving derived LMIs, simulation results are presented and compared with similar earlier researches.

Consider a multi-agent system consisting of three single link manipulators with flexible joints actuated by DC motors, as depicted in Figure 1.

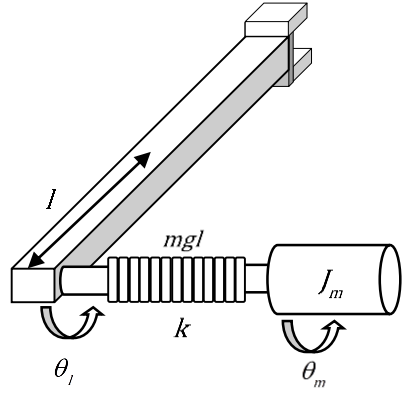

**Figure 1: single link manipulator with flexible joint**

The dynamics of a single link manipulator with flexible joint can be expressed by [23]:

$$J_l \ddot{\theta}_l + mgl \sin(\theta_l) + k(\theta_l - \theta_m) = 0$$
$$J_m \ddot{\theta}_m + k(\theta_m - \theta_l) + b\dot{\theta}_m = k_\tau u \quad (35)$$

where $\theta_m$, $J_m$, $\theta_l$ and $J_l$ represent the motor angular rotation, the motor inertia, the link angular rotation and the link inertia, respectively. Moreover, $m$ is the pointer mass, $l$ is half link length, $k$ is torsional spring constant. Viscous friction coefficient is denoted by $b$ and $k_\tau$ is the amplifier gain. On the other hand, the derivation of angular rotation in links and motors are defined as angular velocity of the form:

$$\omega_l = \dot{\theta}_l, \ \omega_m = \dot{\theta}_m \quad (36)$$

Assuming that only the states of the motors are available, the state space description of a single link manipulator is written according to nominal values in [24]:

$$\dot{x}_i = Ax_i + Bu_i + h_i(t, x_i)$$
$$y_i = Cx_i \quad (37)$$

where $x_i = \left(\theta_m^T \ \omega_m^T \ \theta_l^T \ \omega_l^T\right)^T$ is the state vector and:

$$A = \begin{pmatrix} 0 & 1 & 0 & 0 \\ -48.6 & -1.25 & 48.6 & 0 \\ 0 & 0 & 0 & 1 \\ 19.5 & 0 & -19.5 & 0 \end{pmatrix}, B = \begin{pmatrix} 0 \\ 21.6 \\ 0 \\ 0 \end{pmatrix}, C = \begin{pmatrix} 1 & 0 & 0 & 0 \\ 0 & 1 & 0 & 0 \end{pmatrix}, h_i(t,x_i) = \begin{pmatrix} 0 \\ 0 \\ 0 \\ -3.33\sin(\theta_l) \end{pmatrix}. \quad (38)$$

The selected network topology of the multi-agent system is given by Figure 2.

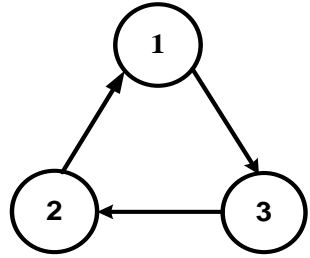

**Figure 2: Network topology of the multi-agent system**

Hence, the Laplacian matrix corresponding to Figure 2 can be written according to (1).

$$L = \begin{bmatrix} 1 & -1 & 0 \\ 0 & 1 & -1 \\ -1 & 0 & 1 \end{bmatrix}. \quad (39)$$

To solve consensus problem using the proposed method in Theorem 1 and corollary 1, in the presence of Gaussian disturbance with zero mean and unit variance, we assume second order controllers with uncertain and certain parameters as two possible scenarios. Let $\hat{H} = I_8$, $R = I_8$ and $Q = 2I_6$ for both investigated methods and bounds $\delta_{A_C} = 0.5$, $\delta_{B_C} = \delta_{C_C} = \delta_{D_C} = 0.2$ in Theorem 1, with the following parameter uncertainties for perturbed dynamic output feedback controllers:

$$\Delta A_C = 0.5 I_3 \otimes \begin{bmatrix} \sin(3t) & \sin(5t) \\ \sin(2t) & \cos(2t) \end{bmatrix}, \quad \Delta B_C = 0.2 I_3 \otimes \begin{bmatrix} \sin(2t) \\ \cos(2t) \end{bmatrix}$$

$$\Delta C_C = 0.2 I_3 \otimes \begin{bmatrix} \cos(t) & \sin(4t) \end{bmatrix}, \quad \Delta D_C = 0.2 I_3 \otimes \sin(t).$$

The obtained controller parameters according to Theorem 1 and Corollary 1, are collected in **Error! Reference source not found.** and the angular rotation and angular velocity of each manipulator is depicted in Figure 3.

**Table 1. Controller parameters for the agents using (I) Theorem 1 and (II) Corollary 1.**

| Parameter | | Agent 1 | Agent 2 | Agent 3 |
|---|---|---|---|---|
| $A_c$ | (I) | $\begin{bmatrix} -1.1021 & 0.0447 \\ 0.0453 & -1.0901 \end{bmatrix}$ | $\begin{bmatrix} -1.0819 & 0.0649 \\ 0.0645 & -1.0628 \end{bmatrix}$ | $\begin{bmatrix} -1.121 & 0.019 \\ 0.019 & -1.122 \end{bmatrix}$ |
| | (II) | $\begin{bmatrix} -1.818 & -0.0097 \\ -0.0099 & -1.859 \end{bmatrix}$ | $\begin{bmatrix} -1.8164 & -0.0087 \\ -0.0087 & -1.8167 \end{bmatrix}$ | $\begin{bmatrix} -1.8184 & -0.0103 \\ -0.0104 & -1.8174 \end{bmatrix}$ |
| $B_c$ | (I) | $\begin{bmatrix} 0.1262 & 0.4236 \\ -1.1873 & -1.1258 \end{bmatrix}$ | $\begin{bmatrix} -0.2333 & -0.1823 \\ 0.3441 & 0.4822 \end{bmatrix}$ | $\begin{bmatrix} -0.0115 & -0.01522 \\ -0.0069 & -0.0027 \end{bmatrix}$ |
| | (II) | $\begin{bmatrix} -0.0066 & -0.0078 \\ -0.0082 & -0.0088 \end{bmatrix}$ | $\begin{bmatrix} -0.0066 & -0.0089 \\ -0.0069 & -0.0038 \end{bmatrix}$ | $\begin{bmatrix} 0.0012 & 0.0044 \\ 0.0015 & 0.0059 \end{bmatrix}$ |
| $C_c$ | (I) | $\begin{bmatrix} -0.0163 & 0.1477 \end{bmatrix}$ | $\begin{bmatrix} 0.0364 & -0.0339 \end{bmatrix}$ | $\begin{bmatrix} 0.0011 & 0.0007 \end{bmatrix}$ |
| | (II) | $\begin{bmatrix} 0.508 & 0.639 \end{bmatrix}*10^{-3}$ | $\begin{bmatrix} 0.627 & 0.641 \end{bmatrix}*10^{-3}$ | $\begin{bmatrix} 0.415 & 0.371 \end{bmatrix}*10^{-3}$ |
| $D_c$ | (I) | $\begin{bmatrix} -0.4987 & -0.5512 \end{bmatrix}$ | $\begin{bmatrix} -0.4951 & -0.5681 \end{bmatrix}$ | $\begin{bmatrix} 0.009 & 0.011 \end{bmatrix}$ |
| | (II) | $\begin{bmatrix} -0.739 & -0.823 \end{bmatrix}$ | $\begin{bmatrix} -0.741 & -0.853 \end{bmatrix}$ | $\begin{bmatrix} 0.003 & 0.009 \end{bmatrix}$ |

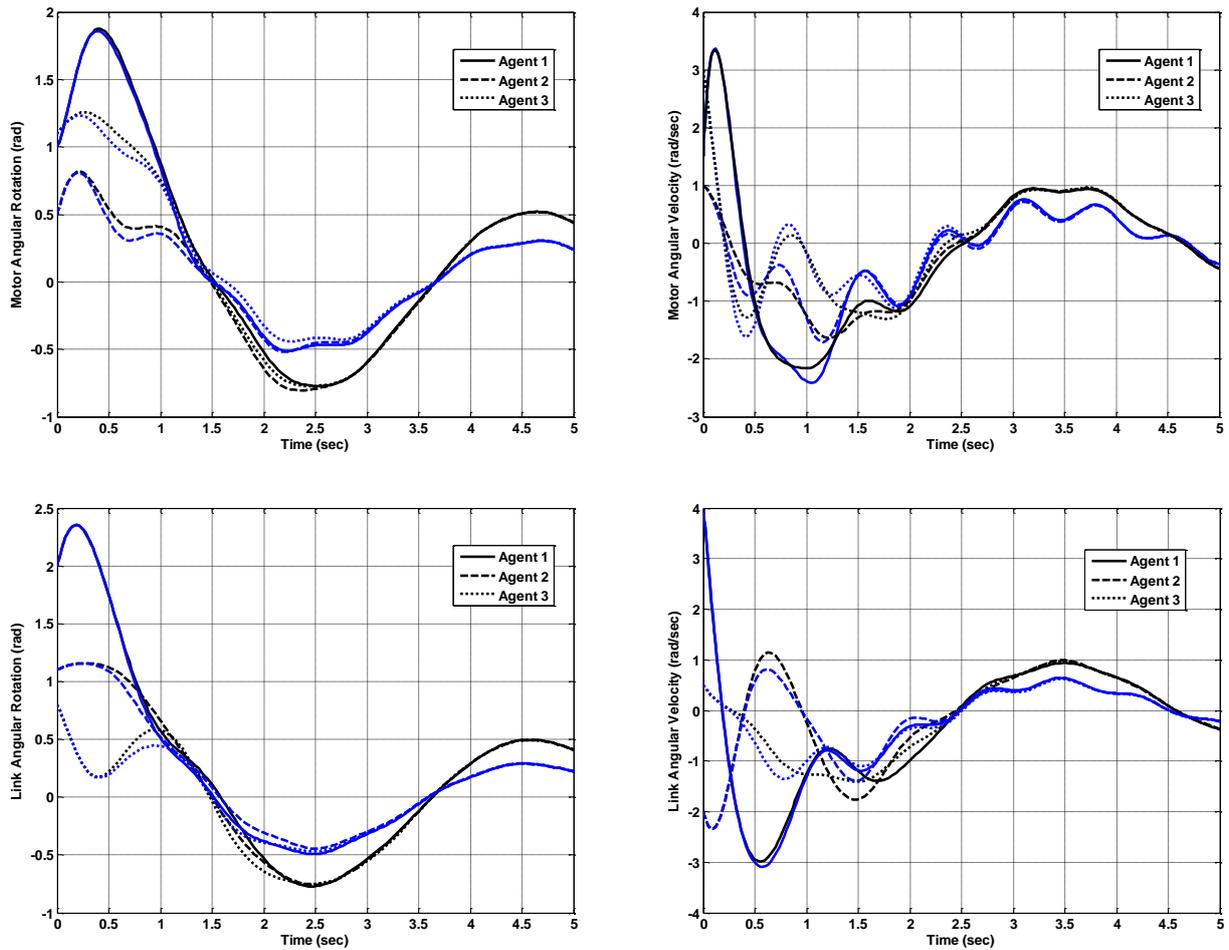

**Figure 3: State consensus in single link manipulator with flexible joint using the proposed controller in Theorem 1 (blue) and Corollary 1 (black).**

According to Figure 3, states convergence of single arm link manipulators has no considerable oscillation using dynamic output feedback controllers. Moreover, consensus is slightly faster when the controllers are assumed to be certain, while the closed-loop multi-agent system response remains reliable with the existence of uncertainties in controllers using Theorem 1.

In order to compare the controller performance with similar works, consensus of Lipschitz nonlinear multi-agent systems with the static state feedback consensus protocol studied in [25] is opted. Simulation results in Figure 4 shows consensus among three single link manipulators using [25] with the same dynamics, communication topology and initial conditions.

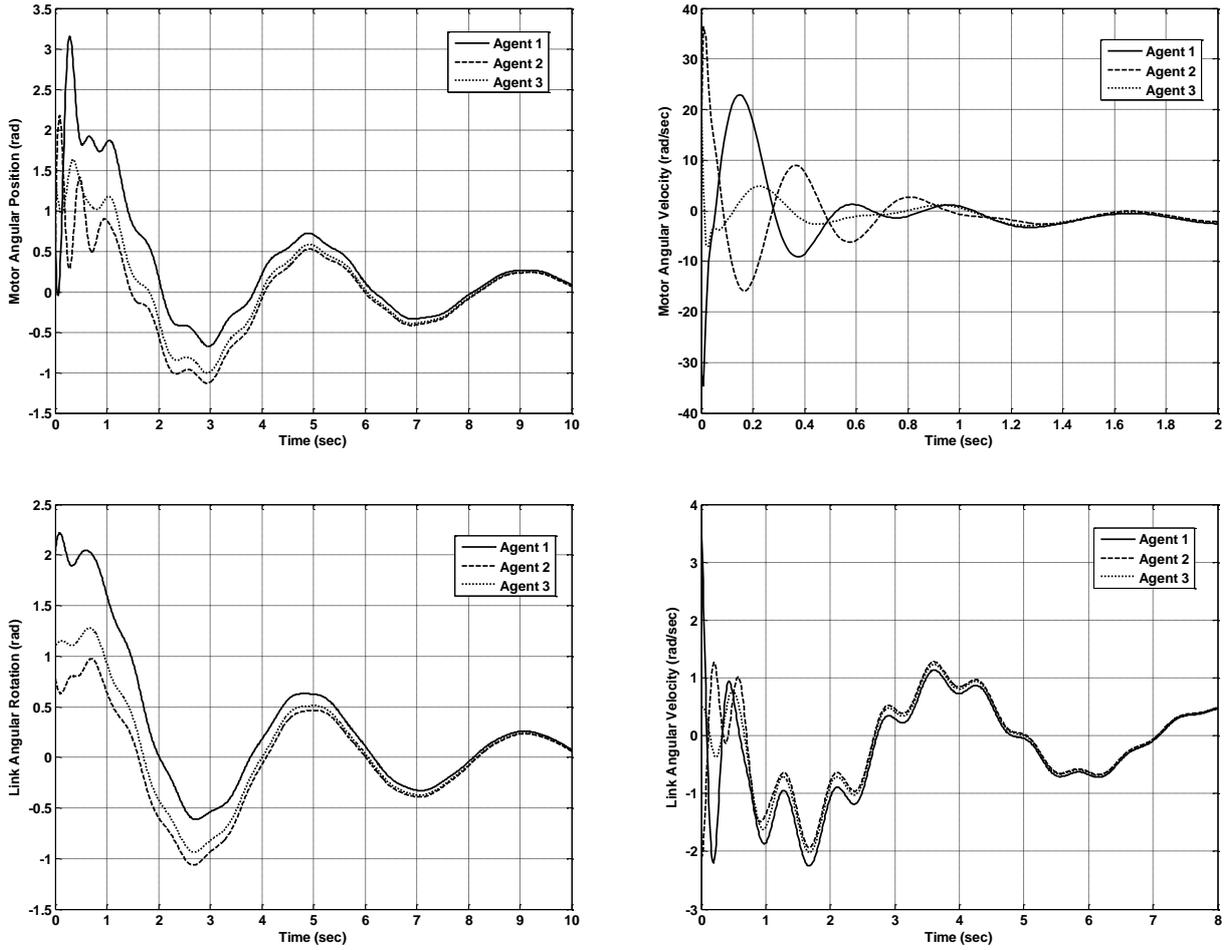

**Figure 4: State consensus in single link manipulator with flexible joint using the proposed controller in [25]**

The trajectories of controller in [25] and the proposed method in Theorem 1, are depicted in Figure 5 because of the importance of control effort and oscillation of actuators in controller implementation.

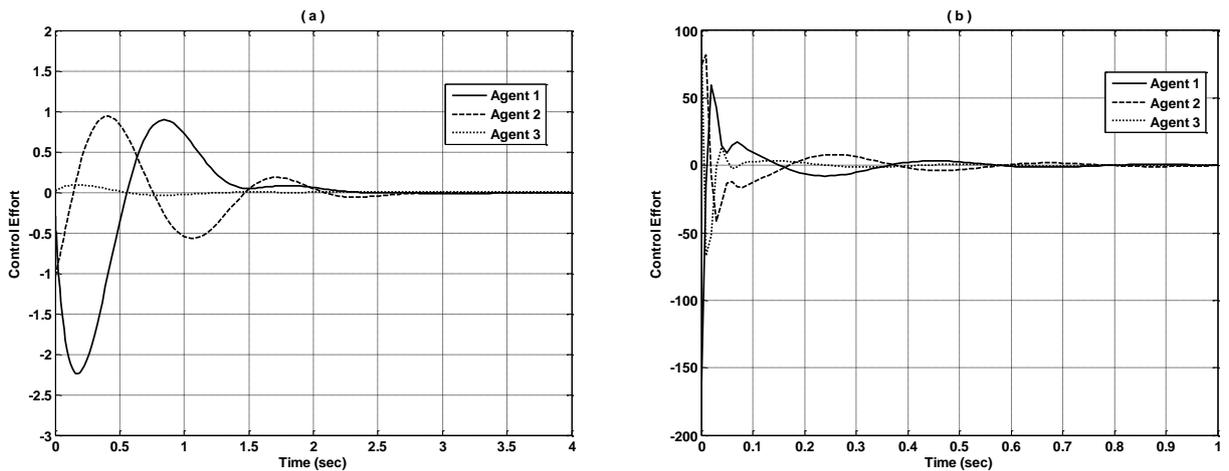

**Figure 5: controller effort of the proposed controller in theorem 1 (a) and [25] (b)**

In addition to a very high overshoot in the beginning instances, the oscillation in controller outputs can be seen in Figure 5 (b), in contrast to a smooth and non-fluctuated controller effort derived by the dynamic output feedback proposed in Theorem 1. Moreover, to analytically compare transient controllers' effort for the proposed protocol in [25] and theorem 1, ISE, IAE, ITSE and ITAE ([26]) are utilized as controller performance indices.

Table 2. Controller performance indices for proposed method in: Theorem 1 and [25].

| Control Method | Agent | ISE | IAE | ITSE | ITAE |
|---|---|---|---|---|---|
| **Theorem 1** | Agent 1 | 1.703 | 1.306 | 0.563 | 0.721 |
| | Agent 2 | 0.457 | 0.807 | 0.259 | 0.644 |
| | Agent 3 | 0.0028 | 0.052 | $8.92*10^{-4}$ | 0.029 |
| **Proposed controller in [25]** | Agent 1 | 334.8 | 5.870 | 4.656 | 1.115 |
| | Agent 2 | 171.4 | 5.628 | 5.538 | 1.210 |
| | Agent 3 | 148.7 | 3.038 | 1.319 | 0.341 |

Perusing the results in Table 2, it is obvious the controller performance by using the proposed method in Theorem 1 is much more acceptable according to all performance indices, even with the presence of controller uncertainties and external disturbance.

## 5. Conclusion

In this paper a consensus method for a group of nonlinear multi-agent system using fixed-order non-fragile dynamic output feedback controller is proposed. We reduced the consensus problem to an equivalent multi-agent system which is supposed to be stabilized and an appropriate Lyapunov function is adopted to derive unknown decentralized controllers' parameters, via an LMI optimization. The robust stability of closed loop augmented system is achieved with the presence of nonlinear dynamics of the agents. Moreover, the controller design procedure is formulated to tolerate norm bounded perturbations in parameters and a $H_\infty$ performance index is utilized to attenuate external disturbance. Finally, to demonstrate the effectiveness of the proposed algorithm and compare with similar earlier researches, a numerical example on a multi-agent system consisting of single link flexible manipulators is carried out.